\documentclass[reqno,12pt]{amsart}
\usepackage{amsfonts}
\usepackage[latin1]{inputenc} 
\usepackage[dvips]{color}

\scrollmode

\newtheorem{theorem}{Theorem}[section]
\newtheorem{proposition}[theorem]{Proposition}
\newtheorem{lemma}[theorem]{Lemma}
\newtheorem{corollary}[theorem]{Corollary}

\newcommand{\nn}{\nonumber}

\newcommand{\ra}{\rightarrow}

\newcommand{\Z}{\mathbb{Z}}

\makeatletter
\@addtoreset{equation}{section}
\makeatother

\renewcommand\thefigure{\thesection.\@arabic\c@figure}
\renewcommand\thetable{\thesection.\@arabic\c@table}

\bibdata{prob}
\bibstyle{alpha}

\title[CLT for Self-repelling Random walks]{Central Limit Theorem for the self-repelling random walk with directed edges}

\author{T. Mountford, G. Valle, L. P. R. Pimentel}

\date{\today}

\address{
\newline
Thomas Mountford
\newline
D\'epartement de Math\'ematiques, \'Ecole Polytechnique F\'ed\'erale,
\newline 1015 Lausanne, Switzerland.
\newline
e-mail: {\rm \texttt{thomas.mountford@epfl.ch}}
\newline
\newline
Leandro P. R. Pimentel
and Glauco Valle
\newline
UFRJ - Departamento de m\'etodos estat\'{\i}sticos do Instituto de Matem\'atica,
\newline  Caixa Postal 68530, 21945-970, Rio de Janeiro, Brasil.
\newline
e-mail: {\rm \texttt{leandro@im.ufrj.br}}
\newline
e-mail: {\rm \texttt{glauco.valle@im.ufrj.br}}
}


\begin{document}

\maketitle

\begin{abstract}
We prove a conjecture of Toth and Veto \cite{tv} about the weak convergence of the self repelling random walk with directed edges under diffusive scaling to a uniform distribution. 
\end{abstract}

\section{Introduction}

In their paper, \cite{tv}, Toth and Veto study a self repelling random walk on  $\mathbb{Z}$. To define this process fix a non-decreasing function $w: \mathbb{Z} \rightarrow \mathbb{R}_+$ such that $\lim_{z \rightarrow \infty} \big( w(z) - w(-z) \big) > 0$. The Self-repelling random walk (SRRW) associated to $w$ is a nearest-neighbor random walk $(X(k))_{k \ge 0}$ starting at $X(0) = 0$ and evolving according to the following transition probabilities: 
\begin{eqnarray}
\label{definition}
\lefteqn{
P\big( \, X(k+1) = X(k) \pm 1 \, \big| X(0),...,X(k) \big) \, = \, } \nonumber \\
& & \qquad \frac{ w\big( \mp \big( l^+(k,X(k)) - l^-(k,X(k)) \big) \big) }{ w\big( l^+(k,X(k)) - l^-(k,X(k)) \big) + w\big( l^-(k,X(k)) - l^+(k,X(k)) \big) } \, .
\end{eqnarray}
where $l^+ (k,x)$ and $l^- (k,x)$ are respectively the local times of the directed edges $x \rightarrow x+1$ and $x \rightarrow x-1$, i.e.
$$
l^\pm (k,x) = \# \big\{ 0 \le j \le k-1 : \, X(j)=x, \ X(j+1) = x \pm 1 \big\} \, ,
$$

Many results were proven for this process which was analyzed in \cite{tv} by clever Ray-Knight arguments
following the blueprint of \cite{t}. In particular, it was shown in \cite{tv} that 
$$
\lim_{k \rightarrow \infty} \sup_x \left| \frac{l^+(k,x)}{\sqrt{k}} - \frac{1}{2}\Big( 1 - \frac{x}{\sqrt{k}} \Big)_+ \right| = 0 \, .
$$
in probability, for precise statements see Theorem 1 and corollary 1 on \cite{tv}.
Given this result, it is natural to conjecture that $X(k) / \sqrt k$ converges in distribution to the uniform 
distribution on $(-1,1)$. In \cite{tv} it is also shown that this is the only possible 
nontrivial limit after renormalization. The purpose of this paper is to prove this conjecture and our main result
is the following:

\begin{theorem} \label{theorem1}
Let $(X(k))_{k \ge 0}$ be the SRRW as described above. We have that as $k  \rightarrow \infty$,
$\frac{X(k)}{\sqrt k} $ converges in distribution to $U(-1,1)$, the uniform distribution on $(-1,1)$.
\end{theorem}

{\it Remark:} While this paper was being written we discovered the recent article of \cite{dt} 
in which the analogous result is shown for the self repelling random walk with undirected edges. It is also worth mention that similar questions arise for random walks with site repulsion, see \cite{tv2}.

\section{Notation and preliminary results}
\label{sec:notation}

The paper \cite{tv} is the main reference in this work and, as will be clear throughout the paper, we rely constantly on the ideas and techniques found there. Thus, in order to aid the reader, we will keep notation as close as possible to those in \cite{tv}.

\bigskip

The definition in (\ref{definition}) of the SRRW leads us naturally to a Ray-Knight approach in order to obtain results for the SRRW. This was exploited in \cite{tv}. The main tool is a representation of the local times on the inverse local times,
$$
T^{\pm}_{x,m} = \min \big\{ k\ge 0 : l^{\pm}(k,x) = m \big\} \, ,
$$ 
in terms of independent ergodic Markov chains. We now describe this representation: For every $x \in \mathbb{Z}$ denote
\begin{eqnarray*}
\gamma_x (m) & := & \min\{ k \ge 0 : l^+(k,x) + l^-(k,x) \ge m \} \\
& = & \left\{
\begin{array}{ll} 0 & , \quad m=0 \\ 1 + \min\{ k \ge 0  : \#\{0 \le j \le k : X(j) = x \} = m \}  &, \quad m \ge 1 \, , \end{array}
 \right.
\end{eqnarray*}
$\tau^{\pm,x}_0 = 0$ and for $j \ge 1$
\begin{eqnarray*}
\tau^{\pm,x}_j & = &
\min \big\{ m > \tau^{\pm,x}_{j-1} : l^+\big( \gamma_x(m) , x \big) - l^-\big( \gamma_x(m) , x \big) = \\
& & \qquad \qquad \qquad \qquad \qquad l^+\big( \gamma_x(m-1) , x \big) - l^-\big( \gamma_x(m-1) , x \big) \pm 1 \big\} \, .
\end{eqnarray*}
Put
$$
\eta^{+,x}_j = - \Big( l^+\big( \gamma_x(\tau^{+,x}_j) , x \big) - l^-\big( \gamma_x(\tau^{+,x}_j) , x \big) \Big) \, ,
$$
$$
\eta^{-,x}_j = l^+\big( \gamma_x(\tau^{-,x}_j) , x \big) - l^-\big( \gamma_x(\tau^{-,x}_j) , x \big)  \, ,
$$
for every $j \ge 0$. The Proposition 1 in \cite{tv} states that the processes $\big( \eta^{+,x}_j \big)_{j \ge 0}$, $x\in \mathbb{Z}$, are iid Markov chains starting with initial condition $\eta^{+,x}_0 = 0$. Moreover, from Lemma 1 also in \cite{tv}, these Markov chains are exponentially ergodic, meaning exponencial convergence to the stationary distribution, and their stationary distribution $\nu = \nu_w$ can be described explicitly in terms of $w$, see \cite{tv}. The stationary distributions $\nu_w$ have an invariant property with respect to $w$, their mean are equal to $-1/2$. To simplify notation, we will denote
$$
r^{+,x}_j = \eta^{+,x}_j + \frac{1}{2} \qquad \textrm{and} \qquad r^{-,x}_j = \eta^{-,x}_j + \frac{1}{2} \, .
$$
We will denote by $(r_x)_{x\in \mathbb{Z}}$ a family of iid random variables distributed as $\eta + \frac{1}{2}$ with $\eta$ having distribution $\nu$. This common distribution of the $r_x$'s have the following properties: 
\begin{enumerate}
\item[(i)] the distribution is symmetric, in particular the mean is zero;
\item[(ii)] there are some exponential moments;
\item[(iii)] all integer values have positive mass and the distribution is aperiodic.
\end{enumerate}
We denote the variance of these random variables by $\sigma^2$. 

Below we summarize some useful relations between the local times, inverse local times and the the processes $\big( \eta^{+,x}_j \big)_{j \ge 0}$. For $x\le 0$ and $m \in \mathbb{N}$, we have
\begin{equation}
\label{l1}
l^+(T^+_{x,m},x) = m \, ,
\end{equation}
\begin{equation}
\label{l2}
l^+(T^+_{x,m},y+1) = l^+(T^+_{x,m},y) + 1 + \eta^{-,y+1}_{l^+(T^+_{x,m},y) + 1} \quad x \le y < 0 \, , 
\end{equation}
\begin{equation}
\label{l3}
l^+(T^+_{x,m},y+1) = l^+(T^+_{x,m},y) + \eta^{-,y+1}_{l^+(T^+_{x,m},y)} \quad y \ge 0 \, , 
\end{equation}
\begin{equation}
\label{l4}
l^+(T^+_{x,m},y-1) = l^+(T^+_{x,m},y) + \eta^{+,y}_{l^+(T^+_{x,m},y)} \quad y \le x \, . 
\end{equation}

\bigskip

As a first result, we will improve some bounds in \cite{tv} on the tail of the distributions of the rightmost and leftmost positions visited by the SRRW by times $T^{\pm}_{0,m}$. Let us start with the proper definitions. Let
$$
\rho^{\pm}_m = \sup \big\{ x \in \mathbb{Z} : l^+(T^{\pm}_{0,m},x) > 0 \big\} 
$$
and
$$
\lambda^\pm_m = \inf \big\{ x \in \mathbb{Z} : l^-(T^{\pm}_{0,m},x) > 0 \big\} \, .
$$
See Theorem 1 in \cite{tv} for the scaling limit for these quantities.

To simplify notation, put
$$
l^+_m (x) =  l^+( T^{+}_{0,m} , x ) \, .
$$
From a close examination of the proofs of (46), (49), (50), (51) and (52) in \cite{tv} (see details below), we have that for any function $g: \mathbb{N} \rightarrow \mathbb{R}_+$ such that
$\lim_{m \rightarrow \infty} g(m) = + \infty$ there exists constants $\beta > 0$ and $c = c(g) > 0$ such that
\begin{equation}
\label{l>}
P \big( \, l^+_m \big( 2m - 4 \sqrt{m g(m)} \big) \ge  3 \sqrt{m g(m)} \big) \le c \, e^{- \beta g(m)} \, ,
\end{equation}
\begin{equation}
\label{l<}
P \Big( \, \min_{1 \le x \le 2m - 4 \sqrt{m g(m)} } l^+_m(x) \le  \sqrt{m g(m)} \, \Big) \le c \, m \, e^{- \beta g(m)} \, ,
\end{equation}
$$
P \big( \, \rho^+_m \ge 2m + \sqrt{m} g(m) \, \big) \le \frac{c}{\sqrt{g(m)}} 
$$
and
$$
P \big( \, \lambda^+_m \le -( 2m + \sqrt{m} g(m)) \, \big) \le \frac{c}{\sqrt{g(m)}} \, .
$$
The last two inequalities will be proved below. Inequality (\ref{l>}) follows from (50) in \cite{tv}, indeed start replacing $A^{\epsilon}$ by $\sqrt{g(A)}$ exclusively where $\epsilon$ appears, after that take $x=0$, $h=1$ and $A=m$. To obtain inequality (\ref{l<}) combine (46) and (49) in \cite{tv} considering the same especifications for $x$, $h$ and $A$ as above.

\bigskip

We start improving the last two bounds on the tail of the distribution of $\rho^+_m$ and $\lambda^+_m$.

\bigskip

\begin{lemma} 
\label{rho}
For any function $g: \mathbb{N} \rightarrow \mathbb{R}_+$ such that
$\lim_{m \rightarrow \infty} g(m) = + \infty$ and $\limsup_{m \rightarrow \infty} \frac{g(m)}{\sqrt{m}} = 0$ there exist constants $\beta> 0$ and $c = c(g) > 0$ such that
$$
P \big( \, \rho^+_m \ge 2m + \sqrt{m} g(m) \, \big) \le c \, e^{- \beta g(m)} \, .
$$
\end{lemma}

\bigskip

Since Lemma \ref{rho} is not our main concern, although we use it ahead, we postpone the proof to Appendix \ref{ap:rho}. 

\bigskip \bigskip


\section{The local central limit Theorem}

\medskip

In the following we consider the position of the random walk $(X(k))_{k \geq 0}$ at some large time.
In order to better conceive the quantities involved, we write this time as $n^2$ even though, obviously 
a typical positive integer is not a perfect square. The arguments presented will not make use of this
and for a general integer time $k$ the term $n$ should be thought of as the integer part of $\sqrt{k}$. 

A second hypotheses concerns the parity of $x$. Since our aim is to estimate $P\big( X_{n^2} = x \big)$ and
the random walk $(X(k))_{k \geq 0}$ has period 2, the parity of $x$ will have an effect on our expressions. We
will suppose that $x$ and $n^2$ are even and the case $x$ and $n^2$ odd can be dealt analogously.

With the above conventions, Theorem \ref{theorem1} is a straightforward consequence of the following local central limit theorem
for the self-repelling random walk $(X(k))_{k \geq 0}$:

\begin{proposition}
\label{lcltX}
There exists $1/2 < \alpha < 1$ such that, for every $\epsilon > 0$, we can take $n_0 = n_0(\epsilon)$ sufficiently large so that if $n \ge n_0$ then
$$
P(X_{n^2} = x)  \geq  \frac{1-  \epsilon}{n} \, ,
$$ 
for every $|x| \leq n- n^ \alpha$ with the same parity as $n^2$.
\end{proposition}

The rest of this section is devoted to the proof of Proposition \ref{lcltX}. By symmetry of the self-repelling walk, we can consider only the case $x \le 0$ and we will suppose this to be the case. 

\smallskip

Note that $P\big( X_{n^2} = x \big)$ is equal to
\begin{eqnarray}
\label{reprlt}
\lefteqn{P \big( \exists \, 0 \le m \le n^2 \textrm{ such that } T^+_{x-1,m} = n^2 \big) + } \nonumber \\
& & \qquad  P \big( \exists \, 0 \le m \le n^2 \textrm{ such that } T^-_{x+1,m} = n^2 \big) \nonumber \\
& & \qquad \qquad = \sum_{m=0}^{n^2} P(T^+_{x-1,m} = n^2) + \sum_{m=0}^{n^2} P(T^-_{x+1,m} = n^2) \, . 
\end{eqnarray}
Our first step is to consider for which values of $m$ the contribution of $P(T^+_{x-1,m} = n^2)$ is relevant in the sum above. We claim that $m$ should be $n/2$ plus a term of order $\sqrt{n}$, otherwise the contribution of $P(T^+_{x-1,m} = n^2)$ can be neglected. Indeed this is the content of the Lemma \ref{there} which also aims at providing precise asymptotics for $P(T^+_{x-1,m} = n^2)$ for the right $m$. Before we state the result we need to fix some notation. Recall that $\sigma^2$ is the variance of the stationary distribution $\nu$. Also define 
$$
\theta_u (v) := \frac{u}{2} \Big( 1 - \frac{|v|}{u} \Big) \, , \ u > 0 \, , \ v \in \mathbb{R} \, .
$$

\medskip

\begin{lemma} \label{there}
There exists $1/2 < \alpha < 1$ such that,
for every $ \varepsilon > 0$ and $K>0$, there exists $n_0 = n_0(\epsilon,K)$ sufficiently large such that 
$$
 \sqrt{\beta_n \pi } \, n^{3/2} P \big(T^{\pm}_{x,\theta_n (x) + c\sqrt{n} } = n^2 \big) \, \ge \,  
e^{- \frac{4 \, c^2}{\beta_n} } \, - \,  \varepsilon \, ,
$$
for all $n \ge n_0$, $|x| \leq n- n^ \alpha$ with the same parity as $n^2$ and $c \in \big\{ \tilde{c} \in (-K,K) : \theta_n (x) + \tilde{c} \sqrt{n} \in \mathbb{N} \big\} $, where
$$
\beta_n = \frac{2\sigma^2 \big( (1+\frac{|x|}{n})^3+(1-\frac{|x|}{n})^3  \big)}{3} \, .
$$
\end{lemma}

\bigskip

We postpone the proof of Lemma \ref{there} to Section \ref{sec:there} and now we show how the lemma is used to establish Proposition \ref{lcltX}. 

\bigskip

\noindent \textbf{Proof of Proposition \ref{lcltX}:} 
Put 
$$
m_K^1 = \min \Big\{ m : \frac{m - \theta_n (x)}{\sqrt{n}} > -K \Big\} \quad \textrm{and} \quad m_K^2 = \max \Big\{ m : \frac{m - \theta_n (x)}{\sqrt{n}} < K \Big\} \, .
$$
Then, in Lemma \ref{there}, we can write $c$ as $\frac{m - \theta_n (x)}{\sqrt{n}}$ for some $m_K^1 \le m \le m_K^2$. Therefore, If the inequality in the same lemma holds, we have that
\begin{eqnarray*}
\sum_{m = m_K^1}^{m_K^2} P(T^\pm_{x-1,m} = n^2) & \ge & \sum_{m = m_K^1}^{m_K^2} \Big( \frac{ e^{- \frac{4}{\beta_n} \, \frac{(m-\theta_n(x))^2}{n} } }{ \sqrt{\beta_n \pi } \,  n^{3/2}} - \frac{ \varepsilon }{\sqrt{\beta_n \pi } \, n^{3/2}} \Big) \\
& = & \frac{1}{2 \, n} \sum_{j = -K \sqrt{n}}^{K \sqrt{n}} \sqrt{ \frac{ 4 }{\beta_n \pi n } } \ e^{- \frac{4}{\beta_n}  \big( \frac{j^2}{n} \big) }  
\ - \ \frac{ \varepsilon \, K }{\sqrt{\beta_n \pi }} \, \frac{1}{n}  \, .
\end{eqnarray*}
Since 
$$
\lim_{K \rightarrow \infty} \lim_{n \rightarrow \infty} \sum_{j = -K \sqrt{n}}^{K \sqrt{n}} \sqrt{ \frac{ 4 }{\beta_n \pi n } } \ e^{- \frac{4}{\beta_n} \big( \frac{j^2}{n} \big) } = \int_{-\infty}^{\infty} \sqrt{ \frac{ 3 }{\sigma^2 \pi } } e^{- \frac{3 u^2}{\sigma^2}}  du  = 1 \, ,
$$
for each $\epsilon > 0$ fixed, we can choose $K = K(\epsilon)$ sufficiently large and then $\varepsilon = \varepsilon(K,\epsilon)$ such that  
$$
\sum_{m = m_K^1}^{m_K^2} P(T^\pm_{x-1,m} = n^2) \ge \frac{1-\epsilon}{2n} \, ,
$$
for all $|x| \leq n- n^ \alpha$ with $n$ sufficiently large.

By (\ref{reprlt}), we get 
$$
P\big( X_{n^2} = x \big) \ge \frac{1-\epsilon}{n} \, .
$$
$\square$

\bigskip

\section{Proof of Lemma \ref{there}}
\label{sec:there}

\medskip

Recall that we are supposing $x \le 0$ and that $x$ and $n^2$ have the same parity. Express $T^+_{x,m}$ in terms of the ``onward" local times $l^\pm (k,x)$ as 
$$
T^+_{x,m} = \sum_{y \in \mathbb{Z}} \Big( \, l^+ \big(T^+_{x,m} , y \big) + l^- \big(T^+_{x,m} , y \big) \, \Big) =
2 \sum_{y \in \mathbb{Z}} l^+ \big(T^+_{x,m} , y \big) + |x| \, .
$$
Therefore the probability in Lemma \ref{there} is 
\begin{equation}
\label{eq:suml+}
P \Big( \sum_{y \in \mathbb{Z}} l^+ \big(T^+_{x,\theta_n (x) + c\sqrt{n}} , y \big) = \frac{n^2 - |x|}{2} \Big) \, .
\end{equation}
To deal with this last probability, write
$$\sum_{y \in \mathbb{Z}} l^+ \big(T^+_{x,\theta_n (x) + c\sqrt{n}} , y \big) = Z^n_{x,c} + W^n_{x,c,-} + W^n_{x,c,+} \, ,$$
where 
$$
Z^n_{x,c} \ =  \  \sum_{|y| \leq n - n^{1/2} log(n)}  l^+ \Big(T^+_{x,\theta_n (x) + cn^{1/2} } , y \Big) \, ,
$$
$$
W^n_{x,c,1} \ =  \  \sum_{y > n - n^{1/2} log(n)}  l^+ \Big(T^+_{x,\theta_n (x) + cn^{1/2} } ,y \Big) 
$$
and
$$
W^n_{x,c,2} \ =  \  \sum_{y < - n + n^{1/2} log(n)}  l^+ \Big(T^+_{x,\theta_n (x) + cn^{1/2} } ,y \Big) \, .
$$

\medskip

We start considering $Z^n_{x,c}$. We are going to show that it can be replaced by more convenient random variables that reduces the problem to a local central limit theorem for
sums of independent random variables. We need a proper representation of the local times $l^+ \big(T^+_{x,\theta_n (x) + c\sqrt{n}} , y \big)$ in terms of the processes $\big( \eta^{+,x}_j \big)_{j \ge 0}$, $x\in \mathbb{Z}$.
From equalities (\ref{l1}) to (\ref{l4}), $l^+ \big(T^+_{x,\theta_n (x) + c\sqrt{n}} , y \big) $ is equal to
$$
\left\{
\begin{array}{ll}
\theta_n (x) + c\sqrt{n} + |x| + \sum_{z=x+1}^y \eta^{-,z+1}_{j_{n}(z)}  &, \ \textrm{if} \ 0 \le y \le n - n^{1/2} log (n) \, , \\
\theta_n (x) + c\sqrt{n} + |x-y| + \sum_{z=x+1}^y \eta^{-,z+1}_{j_{n}(z)}  &, \ \textrm{if} \ x \le y < 0 \, , \\
\theta_n (x) + c\sqrt{n} + \sum_{z=y}^{x-1} \eta^{+,z}_{j_{n}(z)}  &, \ \textrm{if} \ - n + n^{1/2} log (n) \le y < x \, ,
\end{array}
\right.
$$ 
that can be rewritten as
\begin{equation}
\label{lformula}
\left\{
\begin{array}{ll}
\theta_n (y) + c\sqrt{n} + \sum_{z=x+1}^y r^{-,z+1}_{j_{n}(z)}  &, \ \textrm{if} \ x \le y \le n - n^{1/2} log (n) \, , \\
\theta_n (y) + c\sqrt{n} + \sum_{z=y}^{x-1} r^{+,z}_{j_{n}(z)}  &, \ \textrm{if} \ - n + n^{1/2} log (n) \le y < x \, ,
\end{array}
\right.
\end{equation}
where
$$
j_n(z) = 
\left\{
\begin{array}{ll}
l^+ \big(T^+_{x,\theta_n (x) + c\sqrt{n}} , \, z \big) &, \ \textrm{if} \  - n + n^{1/2} log (n) \le z < 0  \, , \\
l^+ \big(T^+_{x,\theta_n (x) + c\sqrt{n}} , \, z \big) + 1 &, \ \textrm{if} \ 0 < z \le n - n^{1/2} log (n) \, . \\
\end{array}
\right.
$$

\medskip

From here we need the following:

\medskip

\begin{lemma}
\label{iidreplacement}
There exists a coupling between the processes $\big( r^{+,x}_j \big)_{j \ge 0}$ starting with initial condition $r^{+,x}_0 = 1/2$ and the family $(r_x)_{x\in \mathbb{Z}}$,  introduced in Section \ref{sec:notation}, such that
$$
P \big( \mathcal{A}_n \big) \ge 1 - C  e^{- \beta \, \log^2(n)} \, ,
$$
for some $\beta>0$ and $C>0$, 
where 
$$
\mathcal{A}_n = \Big\{ l^+ \Big(T^+_{x,\theta_n (x) + cn^{1/2} } ,y \Big) = \tilde{l}^+ \Big(T^+_{x,\theta_n (x) + cn^{1/2} } ,y \Big) \, , \, 0 \leq |y|  \leq n - n^{1/2} log(n) \Big\}
$$ 
with
$$
 \tilde{l}^+ \Big(T^+_{x,\theta_n (x) + cn^{1/2} } ,y \Big) = \theta_n (y) + cn^{1/2} + \sum_{z=x+1}^{y} r_z
$$
for $y > x$ and
$$
 \tilde{l}^+ \Big(T^+_{x,\theta_n (x) + cn^{1/2} } ,y \Big) = \theta_n (y) + cn^{1/2} + \sum_{z=y}^{x-1} r_z 
$$
for $y<x$.
\end{lemma}

\medskip

\noindent \textbf{Proof:} The construction of the coupling is described in \cite{tv}. Similarly to (\ref{l<}), by inspection of the proofs of (46) and (49) in \cite{tv}, we have that for any function $g: \mathbb{N} \rightarrow \mathbb{R}_+$ such that
$\lim_{m \rightarrow \infty} g(m) = + \infty$ there exists constants $\beta > 0$ and $c = c(g) > 0$ such that
$$
P \Big( \, \min \Big\{ k : l^+ \big(T^+_{x,j}  ,y \big) = \tilde{l}^+ \big(T^+_{x,j} ,y \big) \, , \forall j \ge n h \Big\} \, \le \, |x| + 2hn - 4 \sqrt{ng(n)} \Big)\, ,
$$
is bounded above by a term of order $n \, e^{- \tilde{\beta} g(n)}$ for some $\tilde{\beta} > 0$. Since $|x| + 2 \theta_n(x) = n $ and $|c| \le M$, choose $h = (\theta_n (x) + cn^{1/2})/n$ and $g(n) = \log^2(n)$ to finish the proof by noting that $|x| + 2hn - 4 \sqrt{ng(n)} < n - n^{1/2} log(n)$ for $n$ sufficiently large. $\square$

\bigskip

Write
$$
V^n_1 (x) \  =  \ \sum_{x < y \leq n - n^{1/2} log(n)} \Big(  \theta_n (y) + cn^{1/2} + \sum_{z=x+1}^{y} r_z \Big)
$$
and 
$$
V^n_2 (x) \  =  \ \sum_{-(n - n^{1/2} log(n)) \leq y < x-1} \Big(  \theta_n (y) + cn^{1/2} + \sum_{z=y}^{x-1} r_z \Big) \, .
$$
By Lemma \ref{iidreplacement}, $Z^n_{x,c} = V^n_1 (x) + V^n_2 (x)$ on $\mathcal{A}_n$.
Using the definition of $\theta_n(y)$, we have that
\begin{eqnarray}
\label{V1+V2}
V^n_1(x) + V^n_2(x) & = & \frac{1}{2} \sum_{ |y| \le n - n^{1/2} log(n)} (n-|y|) - \frac{n-|x|}{2} + 2 c n^{1/2} (n - n^{1/2} log(n)) + \nn \\ 
& & \ + S^n_1(x) + S^n_2(x) \nn \\ 
& = & \frac{n^2}{2} + 2 c \, n^{3/2} \ + \  O(n log^2(n)) + S^n_1(x) + S^n_2(x) \, ,
\end{eqnarray}
where
$$
S^n_1(x) = \sum_{z=x+1}^{n - \sqrt{n} \log(n)} (n- n^{1/2} log(n) - z + 1) \, r_z \, , 
$$
and 
$$
S^n_2(x) = \sum_{ z=-(n- n^{1/2} log(n))}^{x-1}  (n- n^{1/2} log(n)+z +1) \, r_z \ .
$$

\medskip

The next step to estimate the probability in (\ref{eq:suml+}) is to prove that the contribution of $W^n_{x,c,1}$ and $W^n_{x,c,2}$ can be neglected since it is of order smaller than $n^2$.

\medskip

\begin{lemma} 
\label{suml+W}
For $M> 0$ sufficiently large, there exists $\beta> 0$ such that 
$$
\sup_{k=1,2} P \big( W^n_{x,c,k} > M \, n \, log^3(n) \big) \le e^{-\, \beta \, log^2(n)} \, .
$$
\end{lemma}

\medskip

\noindent \textbf{Proof:}
We deal with $W^n_{x,c,1}$. The case of $W^n_{x,c,2}$ is analogous. 

We have, by the anologue of (\ref{l4}) for $x>0$, that $W^n_{x,c,1}$ is equal to 
$$
(\sqrt{n} log(n) + log^2(n) ) L + \sum_{y =  n - n^{1/2} log (n)}^{n + log^2(n)} (y - (n - n^{1/2} log (n))) \, \eta^{+,y}_{l^+ \big(T^+_{n - n^{1/2} log(n), L} , y \big)} \, ,
$$
where $L = l^+ \Big(T^+_{x,\theta_n (x) + cn^{1/2} } ,(n - n^{1/2} log(n)) \Big)$. Thus, given that $L$ is of order smaller or equal to $\sqrt{n} log(n)$ and that $\rho^+ \le n+ log^2 (n)$, which, by (2.5), happens with probability bounded below by $1-e^{-\, \beta \, log^2(n)}$, $W^n_{x,c,1}$ is equal to a term of order $n log^2(n)$ plus
$$
(n^{1/2} log(n) + log^2(n)) \sum_{y =  n - n^{1/2} log (n)}^{n + log^2 (n)} | \eta^{+,y}_{l^+ \big(T^+_{n - n^{1/2} log(n), L} , y \big)} | \, .
$$
Now we proceed as in the proof of Lemma \ref{rho} to show that the previous sum is of order smaller than $\sqrt{n} log^2(n)$ with the required high probability, see (A.2), in appendix A. 
$\square$

\bigskip

We now return to (\ref{eq:suml+}). Similar arguments as the ones above can be applied to the event $\{T_{x,\theta_n (x) + cn^{1/2}}^- = n^2\}$ and we can also replace the sum of the local times of right oriented edges by $V_1^n + V_2^2$. Therefore, By Lemma \ref{suml+W} and Lemma \ref{iidreplacement}, (\ref{eq:suml+}) is the sum of a term that decays as $e^{-\, \beta \, log^2(n)}$ plus
\begin{eqnarray}
\label{eq:suml++}
\lefteqn{
\frac{1}{2} \sum_{j=1}^{M \, n \, log^3(n)} \sum_{j^\prime = 1}^{M \, n \, log^3(n)} 
P \Big( \, W^n_{x,c,1} = j \, , \, W^n_{x,c,2} = j^\prime \Big) } \\
& & \quad P \Big( V^n_1(x) + V^n_2(x)  = \frac{n^2 - |x| - (j+j^\prime)}{2} \, \Big| \, \mathcal{A}_n \, , \, W^n_{x,c,1} = j \, , \, W^n_{x,c,2} = j^\prime \Big) \, . \nonumber
\end{eqnarray}
Define
$$
\tilde{l}_{n,1} = \tilde{l}^+ \Big(T^+_{x,\theta_n (x) + cn^{1/2} } ,n - n^{1/2} log(n) \Big)
$$
and
$$
\tilde{l}_{n,2} = \tilde{l}^+ \Big(T^+_{x,\theta_n (x) + cn^{1/2} } ,-(n - n^{1/2} log(n)) \Big) \, .
$$
Note that $V^n_1(x) + V^n_2(x)$ is conditionally independent of $(W^n_{x,c,1},W^n_{x,c,2})$ given $\tilde{l}_{n,1}$ and $\tilde{l}_{n,2}$ and $\mathcal{A}_n$. From (\ref{lformula}), we have that 
$$
\tilde{l}_{n,1} = \theta_n \big(n - n^{1/2} log(n)\big) + cn^{1/2}  + Y^n_1(x)
$$ 
and
$$
\tilde{l}_{n,2} = \theta_n \big(-(n - n^{1/2} log(n))\big) + cn^{1/2}  + Y^n_2(x)
$$  
where
$$
Y^n_1(x)  \  =  \ \sum_{z=x+1}^{n - n^{1/2} log(n)} r_z \quad \textrm{and} \quad
Y^n_2(x)  \  =  \ \sum_{z=-(n - n^{1/2} log(n))}^{x-1} r_z \ .
$$
Therefore, using (\ref{V1+V2}) and the fact that $W^n_{x,c, \pm}$ are conditionally independent of $S^n_i(x)$ given the $Y^n_i(x)$, we have that (\ref{eq:suml++}) can be rewritten as
\begin{eqnarray}
\label{eq:suml+++}
\lefteqn{\!\!\!\!\!\!\!\!\!\!\!\!\!\!\!
\frac{1}{2} \sum_{j=1}^{M \, n \, log^3(n)} \sum_{j^\prime = 1}^{M \, n \, log^3(n)} \!\!\!\!
P \Big( \, W^n_{x,c,1} = j \, , \, W^n_{x,c,2} = j^\prime \Big) \, 
 \sum_{k} \sum_{k^\prime} \, P \big( Y^n_1(x) = k \, , \, Y^n_2(x) = k^\prime \big) 
 \times } \nonumber \\
& & \times
P \Big( S^n_1(x) + S^n_2(x)  = h_{n}(x) - \frac{(j+j^\prime)}{2} \, \Big| \, Y^n_1(x) = k \, , \, Y^n_2(x) = k^\prime \Big) \, . 
\end{eqnarray}
where $h_{n}(x) = c n^{3/2} -|x|/2 + O(n \log^2 (n))$. Recall from the statement of Lemma \ref{there} that $|c| \le K$ and $|x| \le n$, therefore $h_n(x)$ is of order $n^{3/2}$. We are going to show that conditional probabilities in (\ref{eq:suml+++}) are bounded below by terms of order $n^{3/2}$.

\bigskip

From here, we need local central limit theorems for $(Y^n_1,S^n_1)_{n\ge 1}$ and $(Y^n_2,S^n_2)_{n\ge 1}$.  

\medskip

\begin{lemma} \label{lclt}
For each strictly positive $\varepsilon$, we have that for all
$n$ sufficiently large, uniformly over $0 \le |x| \le  n - n^{1/2} log(n)$ and $(a,b) \in \mathbb{Z}^2$,
\begin{eqnarray*}
\lefteqn{
 \Big|  {2\pi    \frac{\sigma^2}{\sqrt{12}} \, n_x^2 \, P \big( Y^n_1(x)=a, S^n_1(x) =b \big)} } \\
& & \qquad \qquad \qquad \qquad \ - \ exp \Big( -\frac{2}{\sigma^2} \big( \frac{a^2}{n_x}  + \frac{3b^2}{n_x^3} -\frac{3ab}{n_x^2} \big) \Big) \Big| \ < \varepsilon \, ,
\end{eqnarray*}
where $n_x = \big( 1+\frac{|x|}{n} \big) n$.
\end{lemma}

\medskip

The proof of Lemma \ref{lclt} is standard and follows a classical approach as in the proof of the local central limit
theorem for lattice distributions. So we include it in Appendix \ref{ap:lclt} for the sake of completeness.

\bigskip

Lemma \ref{lclt} and the local central limit Theorem for $Y^n_1(x) $ alone yields that

\medskip

\begin{corollary} \label{Sclt}
For each finite $M$ and each strictly positive $\varepsilon$, we have that for all
$n$ sufficiently large, uniformly over $- n + n^{1/2}log(n) \le x \le 0$, $|a|\leq Mn^{1/2}$ and  $|b|\leq Mn^{3/2}$,
\begin{eqnarray*}
\lefteqn{
 \Big|  \frac {\sqrt{2\pi } }{\sqrt{12}}  \sigma n_x^{3/2} P \big( S^n_1(x) =b \, | \, Y^n_1(x)=a \big) } \\ 
& & \qquad \qquad \qquad \qquad - \   
exp \Big( -\frac{6}{\sigma^2} \big(\frac{a}{2n_x^{1/2} }  - \frac{b}{n_x^{3/2}} \big)^2 \Big)   \Big| \ <  \  \varepsilon \, ,
\end{eqnarray*}
where $n_x = \big( 1+\frac{|x|}{n} \big) n$.
\end{corollary}

\medskip

We have a similar result for $(Y^n_2,S^n_2)_{n\ge 1}$ when $0 \leq  x \leq n - n^{\alpha }$. This new constraint on $x$
is to guarantee that $n-n^{1/2}log(n)-x$ is sufficiently large.

\medskip

\begin{corollary} \label{Sclt2}
Let $N_u = n-n^{1/2}log(n)-u$, $u\in \mathbb{R}$. For each finite $M$ and each strictly positive $\varepsilon$, we have that for all
$n$ sufficiently large, uniformly over $- n + n^{\alpha} \le x \le 0$, $|a|\leq MN^{1/2}_x$, $|b|\leq MN^{3/2}_x$,
\begin{eqnarray*}
\lefteqn{
 \Big| \frac {\sqrt{2\pi } }{\sqrt{12}} \sigma N^{3/2}_x P \big( S^n_2(x) =b | Y^n_2(x)=a \big) } \\
& & \qquad \qquad \qquad \qquad - \ 
exp \Big( -\frac{6}{\sigma^2} \big( \frac{a}{2N^{1/2}_x }  - \frac{b}{N^{3/2}_x} \big)^2 \Big) \Big| \ <  \  \varepsilon \, .
\end{eqnarray*}
\end{corollary}

\bigskip

From the previous colloraries we are able to obtain the following:

\bigskip

\begin{lemma}  
\label{cltcond}
For each finite $K$ and for each $\varepsilon > 0$, if $n$ (and therefore $N_x = n-n^{1/2}log(n)-x$) 
is sufficiently large then, whenever $|a| \ \leq \ K n^{1/2} $, $|a^\prime| \ \leq \ K N^{1/2}_x $,
$|b| \ \leq \ \frac{K}{9} n^{3/2} $ and $- n + n^{\alpha} \le x \le 0$,
we have
\begin{eqnarray*}
\lefteqn{
\frac {\sqrt{\pi \beta_n} }{2} n^{3/2} P \Big( S^n_2(x) + S^n_1(x) = b \, \Big| \,Y^n_1(x)=a \, , \, Y^n_2(x)=a^\prime \Big) } \\ 
& & \qquad \qquad \qquad \qquad \ge \ 
exp \Big( -\frac{4}{\beta_n} \big( \frac{a}{2n^{1/2} }+ \frac{a^\prime}{2N^{3/2}_x } - \frac{b}{n^{3/2}} \big)^2 \Big) \, - \, \varepsilon
\end{eqnarray*}
where $\beta_n$ is as in the statement of Lemma \ref{there}.
\end{lemma}

\medskip

We postpone the proof of Lemma \ref{cltcond} to Appendix \ref{ap:cltcond}.

\smallskip

Now we are able to finish the proof of Lemma \ref{there}. By (\ref{eq:suml+++}) and Lemma \ref{cltcond}, we have that 
$$
\sqrt{\beta_n \pi } \, n^{3/2} P \big(T^{\pm}_{x,\theta_n (x) + c\sqrt{n} } = n^2 \big)
$$
is bounded from below by the sum of a negative term that decays as $n^{3/2} e^{-\, \beta \, log^2(n)}$ and
\begin{eqnarray}
\label{eq:suml4}
\lefteqn{\!\!\!\!\!\!\!\!\!\!\!\!\!\!\!
\sum_{j=1}^{M \, n \, log^3(n)} \sum_{j^\prime = 1}^{M \, n \, log^3(n)} \!\!\!\!
P \Big( \, W^n_{x,c,1} = j \, , \, W^n_{x,c,2} = j^\prime \Big) } \nonumber \\
& & \!\!\!\!\!\!\! \sum_{|k| \le K n^{1/2}} \sum_{|k^\prime| \le K N^{1/2}} \, P \big( Y^n_1(x) = k \, , \, Y^n_2(x) = k^\prime \big) 
 \times  \nonumber \\
& & \times \Big[
exp \Big( -\frac{4}{\beta_n} \big( \frac{k}{2n^{1/2} }+ \frac{k^\prime}{2N^{3/2}_x } - \frac{1}{n^{3/2}} \big[ h_{n}(x) - \frac{(j+j^\prime)}{2} \big] \big)^2 \Big) \, - \, \varepsilon \Big] \, . 
\end{eqnarray}

Since 
$$
\frac{1}{n^{3/2}} \big[ h_{n}(x) - \frac{(j+j^\prime)}{2} \big] = c + O\Big( \frac{\log^3(n)}{\sqrt{n}} \Big)
$$
it is straightforward to see that (\ref{eq:suml4}) is bounded below by  $e^{-\frac{4c^2}{\beta_n}} - \varepsilon$ minus a term that goes to zero as $n \rightarrow \infty$. So we have obtained Lemma \ref{there}.

\medskip


\bigskip


\appendix

\section{Proof of Lemma \ref{rho}} \label{ap:rho}

By definition of $\rho^+_m$ and property (\ref{l3})
$$
\rho^+_m = \sup \Big\{ x : m + \sum_{z=1}^x \eta^{-,z}_{l^+_m(z-1)} > 0 \Big\} \, .
$$
So, conditioned to $l^+_m(x_0) = l > 0$, 
$$
\rho^+_m = x_0 + \sup \Big\{ x : l + \sum_{j=1}^k \eta^{-,x_0+z}_{l^+_m(x_0+z-1)} > 0 \Big\}
$$
which has the same distribution as $x_0 + \rho^+_l$. Moreover, It is clear that $\rho^+_m$ is an increasing function of $m$. Now, considering the previous formula with $x_0 = 2m-4\sqrt{mg(m)}$, we have, from (\ref{l<}) for $m$ sufficiently large, that
\begin{eqnarray*}
\lefteqn{
P \big( \, \rho^+_m \ge 2m + \sqrt{m} g(m) \, \big) \le c \, m \, e^{- \beta g(m)} } \\
& & \qquad \qquad \qquad \qquad \qquad  + \, P \big( \, \rho^+_{\sqrt{m g(m)}} \ge \sqrt{m} g(m) + 4 \sqrt{m g(m)} \, \big) \, .
\end{eqnarray*}
We are going to estimate 
$$
P \big( \, \rho^+_{\sqrt{m g(m)}} \ge \sqrt{m} g(m) + 4 \sqrt{m g(m)} \, \big) \, . 
$$
By the exponential convergence to equilibrium of the Markov Chain $\eta$ starting at $0$ (see Lemma 1 in \cite{tv}), we have that $\lim_{n \ra \infty} E[\eta_n] = -\frac{1}{2}$, which is the expected value of the chain in equilibrium. Then, we can fix $l_0 \ge 1$ such that 
$$
\sup_{l \ge l_0} E[\eta_{l}] \le - \frac{1}{4} \, .
$$
Now, put $T := T^{+}_{0,\sqrt{m g(m)}}$ and define the random set
$$
\mathcal{A}_m = \{ 1 \le x \le \sqrt{m} g(m) + 4 \sqrt{m g(m)} \, : \, 0 < l^+(T,x) \le l_0 \} \, .
$$
and $p = \inf_{l \le l_0} P^l(0,-l) > 0$. For each point $x \in \mathcal{A}_m$, independently of any other point in $\mathcal{A}_m$, we have that $\eta^{-,x}_{l^+(T,x)} = - l^+(T,x)$ with probability at least $p$, which implies that $l^+(T,x+1) = 0$, i.e, $l^+(T,z) = 0$ for every $z \ge x$. Therefore,
$$
P\big( \# \mathcal{A}_m \ge 4 \sqrt{m g(m)} \big) \le (1-p)^{4 \sqrt{m g(m)}} \, .
$$
Since $\limsup_{m \rightarrow \infty} \frac{g(m)}{\sqrt{m}} = 0$, the probability in the last expression decays to zero faster than the exponencial of any constant time $g(m)$ as $m$ goes to infinity. Thus, from the last inequality, to prove the statement we only have to care about
\begin{equation}
\label{eq:rhom58A}
P \big( \, \rho^+_{\sqrt{m g(m)}} \ge \sqrt{m} g(m) + 4 \sqrt{m g(m)} \, , \, \# \mathcal{A}_m \le 4 \sqrt{m g(m)} \, \big) \, .
\end{equation}

Also by Lemma 1 in \cite{tv}, we can fix $\alpha > 0$ such that $\max_{1 \le j \le l_0} E[e^{\alpha \, l_0 \, |\eta_j|}]$ is finite. Choose a constant $c > 0$ sufficiently large such that 
$$
\beta^\prime := - c \alpha \, + \, 4 \log \big( l_0 \, \max_{1 \le j \le l_0} E[e^{\alpha \, l_0 \, |\eta_j|}] \big) < 0 \, .
$$
Therefore,
\begin{eqnarray}
\label{eq:sumA}
\lefteqn{
P \Big( \, \sum_{x \in \mathcal{A}_m} \eta^{-,x}_{l^+(T,x-1)} \ge c \, \sqrt{m g(m)} \, , \, \# \mathcal{A}_m \le 4 \, \sqrt{m g(m)} \, \Big) \, \le } \nn \\ 
& & \le \,  P \Big( \, \sum_{x \in \mathcal{A}_m} \sum_{j=1}^{l_0} |\eta^{-,x}_{j}| \ge c \, \sqrt{m g(m)} \, , \, \# \mathcal{A}_m \le 4 \sqrt{m g(m)} \, \Big) \nn \\
& & \le \,  P \Big( \, \sum_{x = 1}^{4 \sqrt{m g(m)} } \sum_{j=1}^{l_0} |\eta^{-,x}_{j}| \ge c \, \sqrt{m g(m)} \, \Big) \le e^{- \beta^\prime \sqrt{m g(m)} } \, ,
\end{eqnarray}
where the last inequality follows from Chebyshev's exponential inequality and the independence of the processes $\eta^{-,x}$, $x \ge 1$.

Returning to the probability in (\ref{eq:rhom58A}), we have that
\begin{eqnarray*}
\lefteqn{
P \big( \,\rho^+_{\sqrt{m g(m)}} \ge \sqrt{m} g(m) + 4 \sqrt{m g(m)} \, , \, \# \mathcal{A}_m \le 4 \sqrt{m g(m)} \big) \, = } \\
& & = \, P \big( \, l^+_{\sqrt{m g(m)}} \big( \sqrt{m} g(m) + 4 \sqrt{m g(m)} \big) > 0 \, , \, \# \mathcal{A}_m \le 4 \, \sqrt{m g(m)} \big) \\
& & = \, P \Big( \,  \sum_{x = 1}^{\sqrt{m} g(m) + 4 \sqrt{m g(m)}} \eta^{-,x}_{l^+(T,x-1)} > 0 \, , \, \# \mathcal{A}_m \le 4 \, \sqrt{m g(m)} \Big) \, .
\end{eqnarray*}
Let $\mathcal{B}_m = \{ 1\le x \le \sqrt{m} g(m) + 4 \sqrt{m g(m)} : x \notin \mathcal{A}_m \}$. By (\ref{eq:sumA}) the last term in the previous expression is bounded above by $e^{- \beta^\prime \sqrt{m g(m)} }$ plus
$$
P \Big( \,  \sum_{x \in \mathcal{B}_m} \eta^{-,x}_{l^+(T,x-1)} > - c \, \sqrt{m g(m)}  \, , \, \# \mathcal{A}_m \le 4 \, \sqrt{m g(m)} \Big) \, .
$$
For $m$ sufficiently large such that $4 \, \sqrt{m g(m)} > 1$, the last probability is less or equal to
\begin{eqnarray}
\label{eq:sumA2}
\lefteqn{
\sum_{j=1}^{4 \, \sqrt{m g(m)}} \, P \big( \, \# \mathcal{A}_m = 4 \sqrt{m g(m)} - j \, \big) \times } \nonumber \\ 
& &  \times P \Big( \,  \sum_{x \in \mathcal{B}_m} \eta^{-,x}_{l^+(T,x-1)} > - c \, \sqrt{m g(m)}  \, \Big| \, \# \mathcal{A}_m = 4 \sqrt{m g(m)} - j \, \Big)  \, .
\end{eqnarray}

\smallskip

We will need the following claim whose proof is a straightforward exercise in probability and is left to the reader:

\medskip

\noindent \textbf{Claim:} If $(W_n)_{n \ge 1}$ is a sequence of independent random variables such that 
\begin{itemize}
\item[(i)] $\sup_{n\ge 1} P(|W_n| = y) \le C \, e^{-c|y|}$, for some $c,C>0$ not depending on $y$;
\item[(ii)] $\sup_{n \ge 1} E[W_n] < 0$.
\end{itemize}
Then for every $\alpha > 0$ sufficiently small
$$
\sup_{n \ge 1} E[ e^{\alpha \, W_n} ] < 1 \, ,
$$
and for every $a_n \rightarrow +\infty$ and $r>0$ there exist $C = C(\alpha,(a_n),r)$ and $\beta^{\prime \prime} = \beta^{\prime \prime} (\alpha,(a_n),r)$ such that  
$$
P \Big( \,  \sum_{k = 1}^{n \, a_n} W_k \ge - r \, n\, \sqrt{a_n}  \, \Big) \le C e^{- \beta^{\prime \prime} \, n \, a_n }
$$

\medskip

Apply the claim above with $n =  \sqrt{m}$, $a_n = 2 g(m)$ and $r = c$ by choosing $\alpha > 0$ appropriately, to show that (\ref{eq:sumA2}) is bounded above by $C e^{- \beta^{\prime \prime}  \sqrt{m} g(m) }$ for some $C,\beta > 0$.

\bigskip


\section{Proof of the Multivariate Local Central Limit Theorem} \label{ap:lclt}

This appendix is devoted to the proof of Lemma \ref{lclt}. It follows the same steps of the 
classical proof of the local central limit theorem for lattice distributions that can be found
for instance in Section 2.5 in \cite{d}.

Let $(\xi_j)_{j\ge 1}$ be a sequence of iid random variables with mean zero, finite variance $\sigma^2$ and distribution concentrated on $\mathbb{Z} + 1/2$. Denote its common characteristic function by $\phi$. We are going to give an idea of the proof that 
\begin{equation}
\label{ltiid}
\lim_{N \rightarrow \infty} \sup_{u,v \in \mathcal{L}_N} \Big| \frac{\pi \sigma^2}{\sqrt{3}} N^2 P \Big( N^{-1/2} \sum_{j=1}^N \xi_j = u \, , \, N^{-3/2}  \sum_{j=1}^N j \, \xi_j = v \Big) - e^{ - \frac{2}{\sigma^2} (u^2 + 3 v^2 + 3 uv)  } \Big| = 0 \, ,
\end{equation}
where $\mathcal{L}_N = \{ (u,v) : u N^{1/2} \in \mathbb{Z} + N/2 \, , \, v N^{3/2} \in \mathbb{Z} + N(N-1)/4 \}$. 
From the inversion formula for the Fourier transform we have that 
$$
e^{ - \frac{2}{\sigma^2} (u^2 + 3 v^2 + 3 uv)  }
\quad
\textrm{ and }
\quad
P \Big( N^{-1/2} \sum_{j=1}^N \xi_j = u \, , \, N^{-3/2}  \sum_{j=1}^N j \, \xi_j = v \Big)
$$
are respectively equal to
$$
\frac{\sigma^2}{4 \sqrt{3} \pi} \int_{-\infty}^{\infty} \int_{-\infty}^{\infty} e^{-itu-isv} \, e^{-\frac{\sigma^2}{2}\big( t^2 + t s + \frac{s^2}{3} \big)} dtds \, ,
$$
and
$$
\frac{1}{4 \pi^2 N^2} \int_{-\pi N^{3/2}}^{\pi N^{3/2}} \int_{-\pi N^{1/2}}^{\pi N^{1/2}} e^{-itu-isv} \phi_N (t,s) dtds \, ,
$$
where $\phi_N$ is the bivariate characteristic function of $\big( N^{-1/2}  \sum_{j=1}^N \xi_j \, , \, N^{-3/2}  \sum_{j=1}^N j \, \xi_j \big)$. Therefore the absolute value expression in (\ref{ltiid}) is bounded above by 
$$
\frac{1}{4 \pi^2 N^2} \int_{-\pi N^{3/2}}^{\pi N^{3/2}} \int_{-\pi N^{1/2}}^{\pi N^{1/2}} \Big| \phi_N (t,s) - e^{-\frac{\sigma^2}{2}\big( t^2 + t s + \frac{s^2}{3} \big)} \Big| dtds  
$$
plus
$$
\frac{\sigma^2}{4 \sqrt{3} \pi} \int \int_{\Gamma_N} e^{-\frac{\sigma^2}{2}\big( t^2 + t s + \frac{s^2}{3} \big)} dtds \, ,
$$
where $\Gamma_N = ([\pi N^{1/2}, \pi N^{1/2}] \times [\pi N^{3/2}, \pi N^{3/2}])^c$. The second integral clearly goes to zero as $N \rightarrow \infty$. The first integral also goes to zero as $N \rightarrow \infty$. To prove this, it is enough to apply the usual analytical methods applied to the study of convergence of characteristic functions related to central limit theorems. Just to point out that we get the proper constants, note that
\begin{eqnarray*}
\log \phi_N (t,s) & = & \sum_{j=1}^N \log \phi \Big( \frac{t}{N^{1/2}} + \frac{j \, s}{N^{3/2}} \Big) \\
& = & - \frac{\sigma^2}{2} \sum_{j=1}^N \Big(  \frac{t}{N^{1/2}} + \frac{j \, s}{N^{3/2}} \Big)^2 + o(N^{-1}) \\
& = & - \frac{\sigma^2}{2} \Big( t^2 + ts + \frac{s^2}{3} \Big) + o(N^{-1}) \, . 
\end{eqnarray*}
We have that (\ref{ltiid}) follows from the claimed convergences.

\smallskip

From (\ref{ltiid}) applied to the sequence of iid random variables $(r_z)_{z\in \mathbb{Z}}$, we get Lemma \ref{lclt} by setting $N = n  + |x|$, $u = a/N^{1/2}$ and $v = b/N^{3/2}$.

\bigskip


\section{Proof of Lemma \ref{cltcond}} \label{ap:cltcond}

\bigskip

\medskip

A priori we will prove a result about convolutions of "Approximate discrete" Gaussians.  Not surprisingly these yield "Approximate discrete" Gaussians.  In the following we consider
$f(x) \ = \ \frac{1}{\sqrt 2 \pi} e^ {-x^2 / 2}$ and $\Phi (x) \ = \ \int_x^{\infty } e^{-u^2 / 2}du$ .

\medskip

\begin{lemma} \label{lemund}
For given $M, \ \epsilon \in (0, \infty )$ with $M \geq 1$ and $\epsilon   \leq   \epsilon_0 < < \ 1$,
there exists $\sigma_0 < \infty$  so that for 
$$
\sigma_0 < \sigma_1 \leq \sigma_2 ,
$$
whenever positive sequences $(\xi_k)_{k \in \mathbb{Z}}$ and  $(\zeta_k)_{k \in \Z}$
satisfy 
$$
\xi_k \  \geq \ \frac{1- \epsilon }{\sigma_2} \,  f\Big(\frac{k}{\sigma_2}\Big)
\mbox{ for } |k|  \leq M \sigma_2
$$
and
$$
\zeta_k \  \geq \ \frac{1- \epsilon }{\sigma_1} \,  f\Big(\frac{k}{\sigma_1}\Big)
\mbox{ for } |k|  \leq M \sigma_1,
$$
then for all $|z| \leq M \sigma_2 / 9$,
$$
\sum_k \xi_{z+k} \zeta_{-k} \ \geq \ (1-\epsilon)^3 \big(1-\Phi(M/4)\big) \frac{1}{\sqrt { \sigma_1 ^2 + \sigma _2 ^2}}
f \Big(\frac{z}{\sqrt {\sigma_1 ^2 + \sigma _2 ^2}}\Big) \, .
$$
\end{lemma}

\noindent \textbf{Proof:}
We first note that uniformly over  $|z| \leq M \sigma_2, |y| \leq M \sigma_1 $,
by the uniform continuity of the function $f$, we have
$$
f\Big( \frac{z+y}{ \sigma_2}\Big) f\Big(\frac{-y}{ \sigma_1}\Big)
\geq \ (1- \epsilon) \int_{y- 1/2}^{y+1/2} f\Big( \frac{z+u}{ \sigma_2}\Big) f\Big(\frac{-u}{ \sigma_1}\Big)du
$$
for $ \sigma_1$ sufficiently large.  
Thus for $ \sigma_1$ large and $|z| \leq M \sigma_2/9 $, we have
$$
\sum_k \xi_{z+k} \zeta_{-k} \ \geq (1- \epsilon)^3 \frac{1}{\sigma_2 \sigma_1}  \int_{- M \sigma_1 /2- 1/2}^{ M \sigma_1/2+1/2} f\Big( \frac{z+u}{ \sigma_2}\Big) f\Big(\frac{-u}{ \sigma_1}\Big)du
$$
which after standard manipulations is equal to
$$
 (1- \epsilon)^3 \frac{1}{\sqrt {\sigma_1 ^2 + \sigma _2 ^2}} f\Big(\frac{z}{\sqrt {\sigma_1 ^2 + \sigma _2 ^2}}\Big) \int_{(- M{\sqrt {\sigma_1 ^2 + \sigma _2 ^2}} /\sigma_2)/2- 1/2}^{M{\sqrt { \sigma_1 ^2 + \sigma _2 ^2}} /\sigma_2)/2+ 1/2} f\Big(u-\frac{z  \sigma_1^2}{\sigma_1^2 +\sigma_2^2}\Big)du
$$
From which the result follows. $\square$

\bigskip

Now, for $n$ and $x$ fixed, put
$$
\xi_k := P \big( S^n_1(x) = k - a n_x | Y^n_1(x) = a \big)
$$
and 
$$
\zeta_k := P \big( S^n_2(x) = k - a^\prime N_x | Y^n_2(x) = a^\prime \big).
$$
By Corollaries \ref{Sclt} and \ref{Sclt2}, we have that the sequences $(\xi_k)_{k \in \mathbb{Z}}$ and  $(\zeta_k)_{k \in \Z}$ satisfy the hypothesis of Lemma \ref{lemund} with
$M =\frac{\sqrt{6}K}{\sigma}$, $\sigma_1 = \frac{\sigma}{\sqrt{12}} \big( 1 + \frac{|x|}{n} \big)^{3/2} n^{3/2}$ and $\sigma_2 = \frac{\sigma}{\sqrt{12}} N_x^{3/2}$. Now apply this lemma to obtain Lemma \ref{cltcond} noting that
$$
\sigma_1^2 + \sigma_2^2 = \frac{\beta_n}{8} + O(n^{5/2} \log(n)) \, .
$$
which is of the same order as $\beta_n \sim O(n^3)$. 

\bigskip \bigskip

\noindent \textit{Acknowledgments:} This work was developed during Thomas Mountford's visit to UFRJ which was partially supported by CNPq Science without Borders grant 402215/2012-5. We would like to thank Maria Eulalia Vares for the efforts to make this visit possible. Both Leandro Pimentel and Glauco Valle were also partially supported by Universal CNPq grant 474233/2012-0 and Glauco Valle was also supported by CNPq grant 304593/2012-5.

\end{document}